\documentstyle[amssymb,12pt]{amsart}

\makeatletter
\@ifundefined{theorem}{\newtheorem{theorem}{Theorem}}{}
\@ifundefined{lemma}{\newtheorem{lemma}[theorem]{Lemma}}{}
\@ifundefined{corollary}{\newtheorem{corollary}[theorem]{Corollary}}{}
\@ifundefined{conjecture}{}{}
\@ifundefined{proposition}{\newtheorem{proposition}[theorem]{Proposition}}{}
\@ifundefined{axiom}{}{}
\@ifundefined{remark}{}{}
\@ifundefined{example}{}{}
\@ifundefined{exercise}{}{}
\@ifundefined{definition}{}{}

\@addtoreset{equation}{section}

\@addtoreset{theorem}{section}

\@addtoreset{example}{section}

\@addtoreset{definition}{section}

\@addtoreset{remark}{section}

\makeatother

\begin{document}
\title[Drinfeld-Sokolov reduction for difference operators]{Drinfeld-Sokolov reduction for difference operators and deformations of $W$%
--algebras. II. General Semisimple Case }
\author{M. A. Semenov-Tian-Shansky}

\address{Universit\'e de Bourgogne, Dijon, France,\\
and Steklov Mathematical Institute, St. Petersburg,}

\author{A. V. Sevostyanov}

\address{Institut of Theoretical Physics, Uppsala University, Uppsala, Sweden, 
\newline and Steklov Mathematical Institute, St.Petersburg}
\date{April 10, 1994}
\ 
\maketitle

\begin{abstract}
The paper is the sequel to \cite{I}. We extend the Drinfeld-Sokolov reduction procedure to q-difference
operators associated with arbitrary semisimple Lie algebras. This leads to a new elliptic deformation of the Lie bialgebra structure on the associated loop algebra. The related classical r-matrix is explicitly described in terms of the Coxeter transformation. We also present a cross-section theorem for q-gauge transformations which generalizes a theorem due to R.Steinberg.
 \end{abstract}

\section*{Introduction}

The present paper is the sequel to \cite{I}; we refer the reader to this
paper for a general introduction. Our goal is to extend the results of \cite
{I} to arbitrary semisimple Lie algebras. As an intermediate step we develop
a group-theoretic framework for abstract difference equations associated
with arbitrary semisimple Lie groups. A similar problem for differential
equations which was solved by Drinfeld and Sokolov is linear, since it
involves only the structure of the corresponding semisimple Lie algebra.
Difference equations lead to the study of specific submanifolds in Lie
groups which are closely related to some of its Bruhat cells. Our main
technical result is a cross-section theorem for the q-gauge transformations
which generalizes a theorem due to R.Steinberg. The reduction scheme
outlined in \cite{I} extends to the abstract case as well and the general
conclusion remains the same: the consistency condition for the reduction
imposes very rigid conditions on the underlying classical r-matrix which fix
it completely. The resulting classical r-matrix is new; it yields an
elliptic deformation of the Lie bialgebra structure on the loop algebras
associated with the so called Drinfeld's new realizations of quantized
affine algebras \cite{nr}, \cite{kh-t}. The explicit characterization of the
set of abstract q-difference operators leads to a very simple formula for
this r-matrix in terms of the Coxeter element of the corresponding Weyl
group. (The role of Coxeter transformations in the theory of q-difference
operators should be compared with the role of the dual Coxeter
transformations in the representation theory of affine Lie algebras at the
critical level which is implicit in \cite{F}, \cite{FR}). The Coxeter
element also plays a key role in the definition and the study of the
generalized Miura transformation.

The structure of the paper is as follows. Section 1 gives a description of
the set of abstract q-difference operators associated to an arbitrary
complex semisimple Lie group $G$. This description is based on the
cross-section theorem referred to above; its proof is postponed until
section 3. In section 2 we define a class of Poisson covariant Poisson
structures on the set of q-difference $G$-valued connections; this
definition is preceded by the description of a class of Lie bialgebra
structures on loop algebras. As compared to \cite{I}, we need more details
on the Poisson theory of q-gauge transformations, including the theory of
the double and the twisted factorization. We then formulate our main theorem
which gives an explicit description of the (unique) classical r-matrix which
is compatible with the Drinfeld-Sokolov reduction for q-difference operators
associated with $G.$ At the end of section 2 we also briefly discuss the
reduction for finite difference operators on the lattice which yields a
definition of the classical lattice $W({\frak g})$-algebra extending the
definition of the lattice Virasoro algebra discussed in \cite{I}. Section 4
contains a description of the generalized Miura transformations. Finally, in
section 5 we compare our formulae with the formulae of Frenkel and
Reshetikhin \cite{FR} for the $sl(n)$ case.

\subsubsection*{Acknowledgements.}

The present paper is part of a joint research project which was started by
E.Frenkel and N.Reshetikhin together with the first author. We are indebted
to B.Kostant who has pointed us the paper \cite{st}; one of the authors
(M.S.T.S) would like to thank G.Arutyunov for useful discussion.

\section{ Abstract q-difference operators}

The following notation will be used throughout the paper. Let $G$ be a
connected simply connected finite-dimensional complex semisimple Lie group, $%
{\frak g}$ its Lie algebra. Fix a Cartan subalgebra ${\frak h}\subset {\frak %
g}\ $and let $\Delta $ be the set of roots of $\left( {\frak g},{\frak h}%
\right) .$ Choose an ordering in the root system; let $\Delta _{+}$ be the
system of positive roots and $\{\alpha _1,...,\alpha _l\},$ $l=rank\,{\frak g%
},$ the corresponding set of simple roots. For $\alpha \in \Delta ,\alpha
=\sum_{i=1}^ln_i\alpha _i$ we define its height by 
\[
\text{height\thinspace }\alpha =\sum_{i=1}^ln_i. 
\]
Let $e_\alpha \in {\frak g}$ be a root vector which corresponds to $\alpha
\in \Delta .$ Let 
\[
{\frak b}={\frak h}+\sum_{\alpha \in \Delta _{+}}{\Bbb C}e_\alpha 
\]
be the corresponding Borel subalgebra and 
\[
\overline{{\frak b}}={\frak h}+\sum_{\alpha \in \Delta _{+}}{\Bbb C}%
e_{-\alpha } 
\]
the opposite Borel subalgebra; let ${\frak n}=[{\frak b},{\frak b}]$ and $%
\overline{{\frak n}}=[\overline{{\frak b}},\overline{{\frak b}}]$ be their
respective nil-radicals. We shall fix a nondegenerate invariant bilinear
form $\left( ,\right) $ on ${\frak g}.$ Let $H=\exp {\frak h},N=\exp {\frak n%
},\overline{N}=\exp \overline{{\frak n}},B=HN,\overline{B}=H\overline{N}$ be
the Cartan subgroup, the maximal unipotent subgroups and the Borel subgroups
of $G$ which correspond to the Lie subalgebras ${\frak h},{\frak n},%
\overline{{\frak n}},{\frak b}$ and $\overline{{\frak b}},$ respectively.

Let $W$ be the Weyl group of $\left( {\frak g},{\frak h}\right) ;$ we shall
denote a representative of $w\in W$ in $G$ by the same letter. Let $%
s_1,...,s_l$ be the reflections which correspond to simple roots; let $s=$ $%
s_1s_2\cdot \cdot \cdot s_l$ be a Coxeter element. Put $N^{\prime }=\{v\in
N;svs^{-1}\in \overline{N}\};$ it is easy to see that $N^{\prime }\subset N$
is an abelian subgroup, $\dim N^s=l.$ Put $M^s=Ns^{-1}N;$ it is clear that $%
Ns^{-1}N=N^{\prime }s^{-1}N$ and that $N^{\prime }\times N\rightarrow
M^s:\left( n^{\prime },n\right) \longmapsto n^{\prime }s^{-1}n$ is a
diffeomorphism.

Let ${\cal G}=LG$ be the loop group of $G;$ the group $G$ will be identified
with the subgroup of constant loops in ${\cal G}.$ Fix $q\in {\Bbb C}$ , $%
\left| q\right| <1,$ and let $\tau $ be the automorphism of ${\cal G}$
defined by $g^\tau (z)=g(qz).$ We shall denote the corresponding
automorphism of the loop algebra $L{\frak g}$ by the same letter.\ Let $%
{\cal C}$ be another copy of ${\cal G}$ equipped with the q-gauge action of $%
{\cal G},$ 
\begin{equation}
{\cal G}\times {\cal C}\rightarrow {\cal C};\left( v,L\right) \longmapsto
v^\tau Lv^{-1}.  \label{gau}
\end{equation}
The space ${\cal C}$ will be sometimes referred to as the space of {\em %
q-difference connections} (with values in $G).$ Let ${\cal M}^s$ be the cell
in ${\cal C}$ consisting of loops with values in $M^s,$ and ${\cal S}={\cal N%
}^{\prime }s^{-1}$

\begin{theorem}
\label{cross}The restriction of the gauge action (\ref{gau}) to ${\cal N}%
\subset {\cal G}$ leaves the cell ${\cal M}^s\subset {\cal C}$ invariant.
The action of ${\cal N}$ on ${\cal M}^s$ is free and \ ${\cal S}$ is a
cross-section of this action.
\end{theorem}

The proof will be given in section \ref{proof}. (Its special case which
corresponds to $G=SL(n)$ is presented in \cite{I}.)

{\em Remark}. A closely related theorem is due to Steinberg \cite{st} who
discussed the $\ $action of an algebraic semisimple Lie group on itself by
conjugations (In other words, the situation considered in \cite{st}
corresponds to the trivial automorphism $\tau =id.)$ Steinberg's theorem
asserts that if $G$ is defined over an algebraically closed field, $%
N^{\prime }s^{-1}\subset G$ is a cross-section of the set of regular
conjugacy classes in \ $G.$ In theorem \ref{cross} we replace the action of
the entire group on itself by the action of its unipotent subgroup on its
affine subvariety ${\cal M}^s.$ In the context of Lie algebras a similar
problem was studied by B.Kostant \cite{K}, again in the case of trivial $%
\tau .\square $

To motivate the above definitions let us discuss the case $G=SL(n).$ Let us
choose an order in the root system of ${\frak sl}(n)$ in such a way that
positive root vectors correspond to lower triangular matrices. We may order
the simple roots and choose the matrices $s_k,k=1,...,n-1,$ representing the
corresponding reflections in such a way that the Coxeter element $%
s=s_1s_2...s_{n-1}$ is represented by the matrix 
\[
s=\left( 
\begin{array}{ccccc}
0 & 0 & \cdot \cdot \cdot & 0 & 1 \\ 
-1 & 0 & \cdot \cdot \cdot & 0 & 0 \\ 
0 & -1 & \cdot \cdot \cdot & 0 & 0 \\ 
\cdot \cdot \cdot & \cdot \cdot \cdot & \cdot \cdot \cdot & \cdot \cdot \cdot
& \cdot \cdot \cdot \\ 
0 & 0 & \cdot \cdot \cdot & 0 & 0 \\ 
0 & 0 & \cdot \cdot \cdot & -1 & 0
\end{array}
\right) . 
\]
Then the group $N^{\prime }$ consists of lower triangular unipotent matrices
of the form 
\[
u=\left( 
\begin{array}{ccccc}
1 & 0 & \cdot \cdot \cdot & 0 & 0 \\ 
0 & 1 & \cdot \cdot \cdot & 0 & 0 \\ 
\cdot \cdot \cdot & \cdot \cdot \cdot & \cdot \cdot \cdot & \cdot \cdot \cdot
& \cdot \cdot \cdot \\ 
0 & 0 & \cdot \cdot \cdot & 1 & 0 \\ 
\ast & * & \cdot \cdot \cdot & * & 1
\end{array}
\right) , 
\]
the set $M^s$ consists of all unimodular matrices of the form 
\[
A=\left( 
\begin{array}{ccccc}
\ast & -1 & 0 & \cdot \cdot \cdot & 0 \\ 
\ast & * & -1 & \cdot \cdot \cdot & 0 \\ 
\cdot \cdot \cdot & \cdot \cdot \cdot & \cdot \cdot \cdot & \cdot \cdot \cdot
& \cdot \cdot \cdot \\ 
\ast & * & * & \cdot \cdot \cdot & -1 \\ 
\ast & * & * & \cdot \cdot \cdot & *
\end{array}
\right) , 
\]
and the set $N^{\prime }s^{-1}$ is the set of all companion matrices of the
form 
\[
L=\left( 
\begin{array}{cccccc}
0 & -1 &  & \cdot \cdot \cdot & 0 & 0 \\ 
0 & 0 & -1 & \cdot \cdot \cdot & 0 & 0 \\ 
0 & 0 & 0 & \cdot \cdot \cdot & 0 & 0 \\ 
\cdot \cdot \cdot & \cdot \cdot \cdot & \cdot \cdot \cdot & \cdot \cdot \cdot
& \cdot \cdot \cdot & \cdot \cdot \cdot \\ 
0 & 0 & 0 & \cdot \cdot \cdot & 0 & -1 \\ 
1 & u_{1} & u_{2} & \cdot \cdot \cdot & u_{n-2} & u_{n-1}
\end{array}
\right) . 
\]
Let us associate to $L$ a first order difference equation 
\[
\tau \cdot \psi +L\psi =0, 
\]
where $\psi =\left( \psi _1,\psi _2,...,\psi _n\right) ^t$ is a column
vector, $\psi _k\in {\Bbb C}\left( \left( z\right) \right) .$ It is easy to
see that its first component $\phi :=\psi _1$ satisfies an n-th order
difference equation, 
\[
\tau ^n\phi +u_{n-1}\tau ^{n-1}\cdot \phi +u_{n-2}\tau ^{n-2}\cdot \phi
+...+u_1\tau \cdot \phi +\phi =0,
\]
and, moreover, $\psi _k=\tau ^{k-1}\phi ,$ $k=1,2,...,n.$ Hence the set $%
{\cal N}^{\prime }s^{-1}$ may be identified with the set ${\cal M}_{n,q}$ of
all n-th order q-difference operators. \footnote{ Note that in this paper we use slightly different
conventions, as compared to \cite{I}; in particular, we consider the
q-gauge action by lower triangular matrices as opposed to upper
triangular matrices} In the general case we set, by definition, ${\cal M}_q(G)={\cal M}%
^s/LN;$ as we shall see, the manifold ${\cal M}_q(G)$ carries a natural
Poisson structure and may be regarded as the spectrum of a classical
q-W-algebra.

{\em Remark.} In \cite{DS} Drinfeld and Sokolov use a similar approach to
define the set of abstract higher order differential operators associated to
a given semisimple Lie algebra; the key observation that motivates their
definition is that for ${\frak g}={\frak sl}(n)$ the matrix 
\[
f=\left( 
\begin{array}{cccccc}
0 & -1 & 0 & \cdot \cdot \cdot & 0 & 0 \\ 
0 & 0 & -1 & \cdot \cdot \cdot & 0 & 0 \\ 
\cdot \cdot \cdot & \cdot \cdot \cdot & \cdot \cdot \cdot & \cdot \cdot \cdot
& \cdot \cdot \cdot & \cdot \cdot \cdot \\ 
0 & 0 & 0 & \cdot \cdot \cdot & 0 & -1 \\ 
0 & 0 & 0 & \cdot \cdot \cdot & 0 & 0
\end{array}
\right) \in {\frak g} 
\]
is a principal nilpotent element. Accordingly, in the context of Lie
algebras the set ${\cal M}^s\subset {\cal G}$ is replaced by the affine
manifold ${\cal M}^f=f+L{\frak b}\subset L{\frak g;}$ the main cross-section
theorem then asserts that the gauge action of the subgroup $LN\subset LG$
leaves ${\cal M}^f$ invariant; its restriction to ${\cal M}^f$ is free and
admits a cross-section which is an affine submanifold in $L{\frak g.}$ This
cross-section provides a model for the space of higher order differential
operators. One may notice that the `near coincidence' of the matrices $%
s^{-1} $ and $f$ is peculiar to the $SL(n)$ case. $\square $

\section{Poisson structures for q-difference equations}

In this section we shall construct a class of Poisson structures on ${\cal G}%
=LG$ and on the space ${\cal C}$ of q-difference connections. Our main
theorem asserts that there is a unique Poisson structure in this class which
is compatible with the Poisson reduction over $LN.$ We start with the
description of a family of Lie bialgebra structures on loop algebras.

\subsection{Factorizable Lie bialgebras associated with $L{\frak g}$}

Let ${\sf g}=L{\frak g}$ be the loop algebra; we equip it with the standard
invariant bilinear form, 
\[
\left\langle X,Y\right\rangle =Res_{z=0}\left( X\left( z\right) ,Y\left(
z\right) \right) dz/z. 
\]
Notice that the automorphism $\tau $ satisfies $\left\langle \tau X,\tau
Y\right\rangle =\left\langle X,Y\right\rangle .$ Let ${\sf d}={\sf g}\oplus 
{\sf g}$ be the direct sum of two copies of ${\sf g}$ with the bilinear form 
\begin{equation}
\left\langle \left\langle \left( X_1,X_2\right) ,\left( Y_1,Y_2\right)
\right\rangle \right\rangle =\left\langle X_1,Y_1\right\rangle -\left\langle
X_2,Y_2\right\rangle .  \label{inner}
\end{equation}
Put ${\sf b}=L{\frak b},\overline{{\sf b}}=L\overline{{\frak b}},{\sf n}=Ln,%
\overline{{\sf n}}=L\overline{\frak n},{\sf h}=L{\frak h}$; let $\pi :{\sf b}%
\rightarrow {\sf b}/{\sf n}$, $\overline{\pi }:\overline{{\sf b}}\rightarrow 
\overline{{\sf b}}/\overline{{\sf n}}$ be the canonical homomorphisms; the
quotient algebras ${\sf b}/{\sf n}$, $\overline{{\sf b}}/\overline{{\sf n}}$
may be canonically identified with ${\sf h.}$ Fix a linear operator $\theta
\in End\,{\sf h}$ which satisfies the following conditions:

\begin{enumerate}
\item  $\left\langle \theta X,\theta Y\right\rangle =\left\langle
X,Y\right\rangle $ for any $X,Y\in {\sf h}$

\item  $I-\theta $ is invertible.
\end{enumerate}

We shall assume, moreover, that $\theta $ extends to an automorphism of $LH$
which we denote by the same letter. Put 
\begin{equation}
{\sf g}_\theta ^{*}=\left\{ \left( X_{+},X_{-}\right) \in {\sf b}\oplus 
\overline{{\sf b}}\subset {\sf d};\overline{\pi }\left( X_{-}\right) =\theta
\circ \pi \left( X_{+}\right) \right\} ;
\end{equation}
let $^\delta {\sf g\subset d}$ be the diagonal subalgebra. The following
assertion is well known.

\begin{proposition}
$\left( {\sf d},^\delta {\sf g},{\sf g}_\theta ^{*}\right) $is a Manin
triple.
\end{proposition}

Thus we get a family of Lie bialgebras $\left( {\sf g},{\sf g}_\theta
^{*}\right) $with common double ${\sf d}={\sf g}\oplus {\sf g}$ parametrized
by $\theta \in End\,{\sf h}$ ${\sf ;}$ all these bialgebras are
factorizable. Let

\[
{\sf g}={\sf n}\dot +{\sf h}\dot +\overline{{\sf n}} 
\]
be the `pointwise Bruhat decomposition' of the loop algebra. Let $%
P_{+},P_0,P_{-}$ be the corresponding projection operators which map ${\sf g}
$ onto ${\sf n},{\sf h},\overline{{\sf n}},$ respectively. The classical
r-matrix associated with $\left( {\sf g},{\sf g}_\theta ^{*}\right) $ is the
kernel of the linear operator $^\theta r_{+}$ $\in Hom\left( {\sf g}_\theta
^{*},{\sf g}\right) ,$%
\begin{equation}
^\theta r_{+}=P_{+}+\left( I-\theta \right) ^{-1}.  \label{rtheta}
\end{equation}

Let us also set $^\theta r_{-}:=-^\theta r_{+}^{*};$ the classical
Yang-Baxter identity implies that both $^\theta r_{+}$ and $^\theta r_{-}$
are Lie algebra homomorphisms from ${\sf g}_\theta ^{*}$ into ${\sf g.}$ In
the definition of the Poisson structures it is sometimes convenient to use
their skew-symmetric combination, 
\begin{equation}
^\theta r=\frac 12\left( ^\theta r_{+}+^\theta r_{-}\right) =\frac 12\left(
P_{+}-P_{-}+\frac{1+\theta }{1-\theta }P_0\right) ;  \label{rskew}
\end{equation}
the `perturbation term' in (\ref{rskew}) is the Cayley transform of $\theta
, $%
\begin{equation}
^\theta r_0=\frac{1+\theta }{1-\theta }P_0.  \label{rnaught}
\end{equation}
The mappings $^\theta r_{\pm }$ give rise to group homomorphisms $^\theta
r_{\pm }:\left( LG\right) ^{*}\rightarrow LG$ (which we denote by the same
letters).

The double ${\sf d}={\sf g}\oplus {\sf g}$ has a natural structure of a
factorizable Lie bialgebra associated with the Manin triple $\left( {\sf d}%
,^\delta {\sf g},{\sf g}_\theta ^{*}\right) $. Hence ${\cal D}={\cal G}%
\times {\cal G}$ is a Poisson Lie group which contains both ${\cal G}$ and
its dual group as Poisson subgroups. More precisely, let $\pi :LB\rightarrow
LB/LN,\pi :L\overline{B}\rightarrow L\overline{B}/L\overline{N}$ be the
canonical projections; the quotient groups $LB/LN,L\overline{B}/L\overline{N}
$ may be canonically identified with $LH$ . Let $^\delta {\cal G}\subset 
{\cal G}\times {\cal G}$ be the diagonal subgroup and ${\cal G}^{*}\subset 
{\cal D}$ the subgroup defined by 
\begin{equation}
{\cal G}^{*}=\left\{ \left( x_{+},x_{-}\right) \in LB\times L\overline{B}%
;\theta \circ \pi (x_{+})=\overline{\pi }(x_{-})\right\}  \label{F}
\end{equation}

\begin{proposition}
(i) $^\delta {\cal G},{\cal G}^{*}\subset {\cal D}$ are Poisson Lie
subgroups with tangent Lie bialgebras $\left( {\sf g},{\sf g}_\theta
^{*}\right) $ and $\left( {\sf g}_\theta ^{*},{\sf g}\right) ,$
respectively. (ii) Almost all elements $x\in {\cal G}$ admit a factorization 
$x=x_{+}x_{-}^{-1}$ where $\left( x_{+},x_{-}\right) \in {\cal G}^{*}.$
\end{proposition}

We shall also need the related notion of {\em twisted factorization}.

\begin{proposition}
\label{twi} Suppose that the automorphism $\tau $ commutes with $\theta .$
Then almost all elements $x\in {\cal G}$ admit a twisted factorization 
\begin{equation}
x=x_{+}^\tau x_{-}^{-1},\;{\rm where\;}\left( x_{+},x_{-}\right) \in {\cal G}%
^{*}.  \label{fac}
\end{equation}
\end{proposition}

Factorizations (\ref{F}) and (\ref{fac}) are unique if we assume that both $%
x $ and $\left( x_{+},x_{-}\right) $ are sufficiently close to the unit
element.

\subsection{Gauge covariant Poisson structures and reduction}

Assume that $\theta \in End\,{\sf h}$ commutes with $\tau .$ In that case
the space ${\cal C}$ of q-difference connections admits a natural Poisson
structure which is covariant with respect to the q-gauge action ${\cal G}%
\times {\cal C}\rightarrow {\cal C}.$ We refer the reader to \cite{RIMS}, 
\cite{double} for its construction which is based on the notion of the
twisted Heisenberg double. For any $\varphi \in C^\infty ({\cal C})$ we
denote by $\nabla \varphi ,\nabla ^{\prime }\varphi $ its left and right
gradients.

\begin{theorem}
(i) For any $\theta \in End\,{\sf h}$ satisfying conditions (1), (2) the
bracket 
\begin{equation}
\begin{array}{r}
\ \left\{ \varphi ,\psi \right\} _\tau =\left\langle ^\theta r\left( \nabla
\varphi \right) ,\nabla \psi \right\rangle +\left\langle ^\theta r\left(
\nabla ^{\prime }\varphi \right) ,\nabla ^{\prime }\psi \right\rangle  \\ 
-\ \left\langle \tau \circ ^\theta r_{+}\left( \nabla ^{\prime }\varphi
\right) ,Y\right\rangle -\left\langle ^\theta r_{-}\circ \tau ^{-1}\left(
\nabla \varphi \right) ,\nabla ^{\prime }\psi \right\rangle ,
\end{array}
\label{taubrack}
\end{equation}
satisfies the Jacobi identity. (ii) Equip the space ${\cal C}$ of
q-difference connections with the Poisson structure (\ref{taubrack}); then
the q-gauge action defines a Poisson mapping ${\cal G}\times {\cal C}%
\rightarrow {\cal C}.$ (iii) The subgroup ${\cal N}=LN\subset {\cal G}$ is
admissible and hence ${\cal N}$-invariant functions form a Poisson
subalgebra in $C^\infty ({\cal C}).$
\end{theorem}

(We refer the reader to \cite{I} for a general definition of admissible
subgroups.) Below we shall need another formula for the Poisson bracket $%
\{,\}_\tau $ which is related to the {\em \ }twisted factorization in ${\cal %
G}.$ Let $\varphi \in C^\infty ({\cal C});$ we define $Z_\varphi \in {\sf g}$
by the following relation: 
\[
r_{+}Z_\varphi -\tau ^{-1}\cdot r_{-}Z_\varphi =\nabla ^{\prime }\varphi . 
\]
Let $h\in {\cal C}$ be an element admitting twisted factorization, $%
h=h_{+}^\tau h_{-}^{-1};$ then 
\begin{equation}
\begin{array}{r}
\left\{ \varphi ,\psi \right\} _\tau (h)=\left\langle Adh_{+}\cdot \tau
^{-1}\cdot Z_\varphi -Adh_{-}Z_\varphi ,\nabla \psi \right\rangle \\ 
-\left\langle \nabla \varphi ,Adh_{+}\cdot \tau ^{-1}\cdot Z_\psi
-Adh_{-}Z_\psi \right\rangle
\end{array}
\label{Pbracket}
\end{equation}

We can now state our second main theorem. For $w\in W$ let $R_w\in End\,L%
{\frak h}$ be the linear operator acting in the loop algebra, 
\[
\left( R_wH\right) \left( z\right) =Ad\,w\cdot \left( H\left( z\right)
\right) . 
\]

\begin{theorem}
\label{main}The quotient space ${\cal M}^s/{\cal N}$ =${\cal S}$ is a
Poisson submanifold in ${\cal C}/{\cal N}$ if and only if the endomorphism $%
\theta $ is given by $\theta =R_s\cdot \tau ,$where $s\in W$ is the Coxeter
element.
\end{theorem}

Notice that the Coxeter transformation satisfies conditions (1), (2) imposed
above; since it preserves the root and weight lattices in ${\frak h},$ it
gives rise to an automorphism of $LH.$ In the case when ${\frak g}={\frak sl}(2)$, theorem \ref{main} was
proved in \cite{I} (Theorem 2); in this case $R_s=-Id.$

For future reference let us write down explicitly the twisted factorization
problem associated with the r-matrix $^\theta r$.

\begin{proposition}
Assume that $\theta =R_s\cdot \tau .$Twisted factorization problem in the
loop group ${\cal G}$ associated with the r-matrix $^\theta r$ amounts to
the relation 
\begin{equation}
x=y_{+}y_{-}^{-1},\;y_{+}\in LB,,y_{-}\in L\overline{B},\;\overline{\pi }%
\left( y_{-}\right) =s\left( \pi \left( y_{+}\right) \right) .  \label{twist}
\end{equation}
\end{proposition}

We shall denote by ${\cal G}^{\prime }\subset {\cal G}$ the open subset
consisting of elements admitting twisted factorization described in (\ref
{twist}).

We shall now explicitly describe the kernel of the corresponding classical
r-matrix. The relevant part of this kernel is the `perturbation term' $r_0$
which was defined in (\ref{rnaught}). In the present situation we have 
\begin{equation}
r_0=\frac{I+R_s\cdot \tau }{I-R_s\cdot \tau }P_0.  \label{Cayley}
\end{equation}
Let $H_p\in {\frak h},p=1,...,l,$ be the eigenvectors of the Coxeter
element, $Ad\,s(H_p)=e^{\frac{2\pi ik_p}h}H_p$ (here $h$ is the Coxeter
number and $k_1,...$ $k_l$ are the exponents of ${\frak g.}$ There exists a
permutation $\sigma $ of the set $\left\{ 1,...l\right\} $ such that the
basis$\left\{ H_{\sigma p}\right\} $ is biorthogonal to $\left\{ H_p\right\}
;$ we may assume that $\left\langle H_p,H_{\sigma p}\right\rangle =1.$Then 
\begin{equation}
\left( r_0X\right) \left( z\right) =\sum_{n=-\infty }^\infty \sum_{p=1}^rz^n%
\frac{1+q^n\exp \frac{2\pi ik_p}h}{1-q^n\exp \frac{2\pi ik_p}h}\left\langle
X_n,H_{\sigma p}\right\rangle H_p,q\in {\Bbb C},\left| q\right| <1.
\label{naught}
\end{equation}
Note that the Lie bialgebra studied by Drinfeld in \cite{nr} (see also \cite
{kh-t}) corresponds to the 'crystalline` limit $q\rightarrow 0$ in (\ref
{naught}); in this case $r_0$ amounts to the Hilbert transform in ${\bf h.}$
It is convenient to write 
\begin{equation}
r_0(q,z)=\sum_{p=1}^r\psi _p(q,z)H_p\otimes H_{\sigma p},  \label{rtensor}
\end{equation}
where 
\begin{equation}
\psi _p(q,z)=\sum_{n=-\infty }^\infty \frac{1+q^n\exp \frac{2\pi ik_p}h}{%
1-q^n\exp \frac{2\pi ik_p}h}z^n.  \label{phip}
\end{equation}
Functions $\psi _p$ satisfy the q-difference equations, 
\begin{equation}
\psi _p\left( q,z\right) +\exp \frac{2\pi ik_p}h\cdot \psi _p\left(
q,qz\right) =\delta \left( z\right) -\exp \frac{2\pi ik_p}h\cdot \delta
\left( qz\right) ,  \label{qdi}
\end{equation}
where 
\begin{equation}
\delta \left( z\right) =\sum_{n=-\infty }^\infty z^n.  \label{dirac}
\end{equation}

\subsection{Proof of Theorem \ref{main}}

We briefly recall the geometric criterion that allows to check that a
submanifold of a quotient Poisson manifold is itself a Poisson manifold. Let 
$M$ be a Poisson manifold, $\pi :M\rightarrow B$ a Poisson submersion.
Hamiltonian vector fields $\xi _\varphi ,\varphi \in \pi ^{*}C^\infty (B),$
generate an integrable distribution ${\frak H}_\pi $ in $TM.$

\begin{proposition}
$\label{poisson}$Let $V\subset M$ be a submanifold; $W=\pi (V)\subset B$ is
a Poisson submanifold if and only if $V$ is an integral manifold of ${\frak H%
}_\pi .$
\end{proposition}

Assume that this condition holds true; let $N_V\subset T^{*}M\mid _V$ be the
conormal bundle of $V;$ clearly, $T^{*}V\simeq T^{*}M\mid _V/N_V.$ Let $%
\varphi ,\psi \in C^\infty (W);$ put $\varphi ^{*}=\pi ^{*}\varphi \mid
_V,\psi ^{*}=\pi ^{*}\psi \mid _V.$ Let $\overline{d\varphi },\overline{%
d\psi }\in T^{*}M\mid _V$ be any representatives of $d\varphi ^{*},d\psi
^{*}\in T^{*}V.$ Let $P_M\in \bigwedge^2Vect\,M$ be the Poisson tensor.

\begin{proposition}
\label{reduce} We have 
\begin{equation}
\pi ^{*}\left\{ \varphi ,\psi \right\} \mid _V=\left\langle P_M,\overline{%
d\varphi }\wedge \overline{d\psi }\right\rangle ;  \label{pbr}
\end{equation}
in particular, the r.h.s. does not depend on the choice of $\overline{%
d\varphi },\overline{d\psi }.$
\end{proposition}

Let us now apply proposition \ref{poisson} in the setting of theorem \ref
{main}. It is sufficient to check that the Hamiltonian vector fields
generated by ${\cal N}$-invariant functions on ${\cal C}$ are tangent to $%
{\cal M}^s$ if and only if $r_0$ is given by (\ref{Cayley}). Let $\varphi
\in C^\infty \left( {\cal C}\right) ^{{\cal N}};$ then $\varphi \left(
v^\tau L\right) =\varphi \left( Lv\right) $ for all $v\in {\cal N},$ and
hence $Z:=\nabla \varphi -\tau \cdot \nabla ^{\prime }\varphi \in {\bf b.}$
Since $\nabla ^{\prime }\varphi (L)=Ad\;L^{-1}\cdot \nabla \varphi (L),$ we
rewrite the Poisson bracket on ${\cal C}$ in the following form: 
\[
\left\{ \varphi ,\psi \right\} _\tau (L)=\left\langle rZ+Z-Ad\;L\cdot r\cdot
\tau ^{-1}\cdot Z+Ad\;L\cdot \tau ^{-1}\cdot Z,\nabla \psi (L)\right\rangle
; 
\]
thus in the left trivialization of $T{\cal C}$ the Hamiltonian field
generated by $\varphi $ has the following form: 
\[
\xi _\varphi \left( L\right) =rZ+Z-Ad\;L\cdot r\cdot \tau ^{-1}\cdot
Z+Ad\;L\cdot \tau ^{-1}\cdot Z. 
\]
Assume that $L\in {\cal M}^s,$ $L=vs^{-1}u,v\in {\cal N}^{\prime },u\in 
{\cal N}.$ Put $Z=Z_0+Z_{+},Z_0\in {\bf h,}Z_{+}\in {\bf n.}$ Then 
\[
\xi _\varphi \left( L\right) =r_0Z_0+Z_0+s^{-1}\tau ^{-1}\cdot
Z_0-s^{-1}\tau ^{-1}\cdot r_0Z_0+X+Ad\;\left( v\cdot s^{-1}\right) \cdot Y 
\]
where $X\in {\bf n}^{\prime },Y\in {\bf n.}$ On the other hand, in the left
trivialization of $T{\cal C}$ the tangent space $T_L{\cal M}^f$ is
identified with ${\bf n}^{\prime }+Ad\;\left( v\cdot s^{-1}\right) \cdot 
{\bf n}$ . Hence $\xi _\varphi $ is tangent to ${\cal M}^f$ if and only its $%
{\bf h}$-component vanishes, i.e., if 
\[
r_0+Id+Ad\,s^{-1}\tau ^{-1}-Ad\,s^{-1}\tau ^{-1}r_0=0, 
\]
which is equivalent to (\ref{Cayley}).

\subsection{Lattice ${\cal W}$-algebras}

Let $\Gamma ={\Bbb Z}/N{\Bbb Z}$ be a finite periodic lattice. Set ${\bf G}%
=G^\Gamma ,{\bf g}={\frak g}^\Gamma .$ Let $\tau $ be the automorphism of $%
{\Bbb G}$ induced by the cyclic permutation on $\Gamma ,x_i^\tau =x_{i+1%
{\rm mod}N}.$ We define lattice gauge transformations by 
\[
g\cdot x=g^\tau xg^{-1}.
\]
The definition of gauge covariant Poisson brackets on the space of
q-difference connections has its obvious lattice counterpart. In an equally
obvious way one may construct a class of Lie bialgebra structures on ${\bf g}%
={\frak g}^\Gamma $ which is compatible with reduction over the unipotent
subgroup ${\bf N}=N^\Gamma .$ Namely, let us consider the `pointwise Bruhat
decomposition' 
\[
{\bf g}={\bf n}\dot +{\bf h}\dot +\overline{{\bf n}}{\bf ,n}={\frak n}%
^\Gamma ,{\bf h}={\frak h}^\Gamma ,\overline{{\bf n}}=\overline{{\frak n}}%
^\Gamma 
\]
and set 
\[
^{{\bf \theta }}{\bf r=P}_{+}-{\bf P}_{-}+\frac{I+{\bf \theta }}{I-{\bf %
\theta }}{\bf P}_0,{\bf \theta }\in End\,{\bf h.}
\]
We shall omit the details and formulate only the lattice counterpart of the
main theorems. Let as usual $s\in G$ be a Coxeter element.

\begin{theorem}
(i) The restriction of the gauge action to ${\bf N}=N^\Gamma $ leaves the
subset ${\bf M}^s={\bf N}s^{-1}N$ invariant. (ii) The restricted action is
free and $S={\bf N}^{\prime }s^{-1},${\bf \ }${\bf N}^{\prime }=N^{\prime
\,\Gamma },$ is its cross-section.
\end{theorem}

We shall assume for simplicity that the lattice length $N$ is relatively
prime with the Coxeter number.

\begin{theorem}
The quotient ${\bf M}^s/{\bf N}$ is a Poisson submanifold in the reduced
space if and only if ${\bf \theta =}R_s\cdot \tau .$
\end{theorem}

{\em Remark.} The condition on $\Gamma $ assures that $I-{\bf \theta }$ is
invertible; it is likely that reduction is possible even without this
assumption, but this question needs further study.$\square $

For $G=SL(2)$ this reduction was studied in detail in \cite{I},
Section 6. It gives a discrete version of the Virasoro algebra, which
is closely connected to the lattice Virasoro algebra of \cite{FV}. 

\section{The Cross-Section Theorem}

\label{proof}

We shall prove the following assertion.

\begin{theorem}
For each \ $L\in {\cal N}^{\prime }s^{-1}{\cal N}$ there exists a unique
element $n\in {\cal N}$ such that $n^\tau \cdot L\cdot s^{-1}\in {\cal N}%
^{\prime }s^{-1}.$
\end{theorem}

Let $C_h\subset W$ be the cyclic subgroup generated by the Coxeter element. $%
C_h$ has exactly $l=rank\,{\frak g}$ different orbits in the root system $%
\Delta ({\frak g},{\frak h}{\frak ).}$ The proof depends on the structure of 
$\ $these orbits; for this reason we have to distinguish several cases.%
\footnote{%
The proofs given below do not apply when ${\frak g}$ is a simple Lie algebra
of type $E_6.$}

\begin{proposition}
The theorem is true for ${\frak g}$ of type $\not A_l.$
\end{proposition}

An elementary proof which is based on the matrix algebra is given in \cite{I}%
; below we present a different proof which uses only the properties of the
corresponding root system.

\begin{lemma}
(i) Each orbit of $C_h$ in $\Delta \left( {\frak g},{\frak h}\right) $%
consists of exactly $h$ elements; one can order these orbits in such a way
that $k$-th orbit contains all positive roots of height $k$ and all negative
roots of height $h-k.$
\end{lemma}

Put 
\[
{\frak n}_k=\bigoplus_{\left\{ \alpha \in \Delta _{+},\;ht\,\alpha
=k\right\} }{\frak n}_\alpha ,N_k=\exp {\frak n}_k,{\cal N}_k=LN_k 
\]
For each $k$ we can choose $\alpha _k\in \Delta _{+}$ in such a way that 
\[
{\frak n}_k=\bigoplus_{p=0}^{h-k-1}{\frak n}_{s^p\left( \alpha _k\right) }. 
\]
Put ${\frak n}_k^p={\frak n}_{s^p\left( \alpha _k\right) },N_k^p=\exp {\frak %
n}_k^p,{\cal N}_k^p=LN_k^p.$

Let $L=v\cdot s^{-1}\cdot u,v\in {\cal N}^{\prime },u\in {\cal N};$ we must
find $n\in {\cal N}$ such that 
\begin{equation}
n^\tau \cdot v\cdot s^{-1}\cdot u=v_0\cdot s^{-1}\cdot n.  \label{nx}
\end{equation}
For any $n\in {\cal N}$ there exists a factorization 
\[
n=n_1n_2\cdot \cdot \cdot n_l,\text{where }n_k\in {\cal N}_k; 
\]
moreover, each $n_k$ may be factorized as 
\[
n_k=n_k^0n_k^1\cdot \cdot \cdot n_k^{h-k-1},\;n_k^p\in {\cal N}_k^p. 
\]

For any $n\in {\cal N}$ the element $n^\tau \cdot v\cdot s^{-1}\cdot u$
admits a representation 
\[
n^\tau \cdot v\cdot s^{-1}\cdot u=\tilde vs^{-1}\tilde u,\;\tilde v\in {\cal %
N}^{\prime },\tilde u\in {\cal N}; 
\]
let 
\[
\tilde u=\overrightarrow{\prod_{k=1}^l}\overrightarrow{\,\prod_{p=0}^{h-k-1}}%
\tilde u_k^p,\;\tilde u_k^p\in {\cal N}_k^p, 
\]
be the corresponding factorization of $\tilde u.$ For $x\in G$ we write $%
s\left( x\right) :=s\cdot x\cdot s^{-1}$ (this notation will be frequently
used in the sequel).

\begin{lemma}
We have $\tilde u_k^p=\tau \cdot s\left( n_k^{p-1}\right) V_k^p,$ where the
factors $V_k^p\in {\cal N}_k^p$ depend only on $u,v$ and on $n_j^q$ with $%
j<k.$
\end{lemma}

Assume now that $n$ satisfies (\ref{nx}); then we have $\tilde v=v_0,\tilde
u=n.$ This leads to the following relations: 
\begin{equation}
\tau \cdot s\left( n_k^{p-1}\right) V_k^p=n_k^p,  \label{recur}
\end{equation}
where we set formally $n_k^{-1}=1.$

\begin{proposition}
The system (\ref{recur}) may be solved recursively starting with $k=1,$ $p=0.
$
\end{proposition}

Clearly, the solution is unique. This concludes the proof for ${\frak g}$ of
type $A_l.$

Let now ${\frak g}$ be a simple Lie algebra of type other than $A_l$ and $%
E_6,$ $l$ its rank.

\begin{lemma}
(i) The Coxeter number $h\left( {\frak g}\right) $ is even. (ii) Each orbit
of $C_h$ in $\Delta \left( {\frak g},{\frak h}\right) $ consists of exactly $%
h$ elements and contains an equal number of positive and negative roots.
(iii) Put 
\[
\Delta _{+}^p=\{\alpha \in \Delta _{+};s^{-p}\cdot \alpha \notin \Delta
_{+}\},\;{\frak n}^p=\bigoplus_{\alpha \in \Delta _{+}^p}{\Bbb C}\cdot
e_\alpha ;
\]
then ${\frak n}^p\subset {\frak n}$ is an abelian subalgebra, $\dim {\frak n}%
^p=l.$
\end{lemma}

When ${\frak g}$ is not of type $D_{2k+1}$ this assertion follows from \cite
{bourb} (Chapter 6, no 1.11, prop. 33 and Chapter 3, no 6.2, corr. 3). For $%
{\frak g}$ of type $D_{2k+1}$ it may be checked directly.

Put $N^p=\exp {\frak n}^p;$ let ${\cal N}^p$ be the corresponding subgroup
of the loop group ${\cal G}.${\em \ }Let $L=v\cdot s^{-1}\cdot u,v\in {\cal N%
}^{\prime },u\in {\cal N};$ we must find $n\in {\cal N}$ such that 
\[
v\cdot s^{-1}\cdot u=n^\tau \cdot v_0\cdot s^{-1}\cdot n^{-1}. 
\]
Put 
\begin{equation}
n=n_1n_2\cdot \cdot \cdot n_{\frac h2},n_p\in {\cal N}_p;  \label{n}
\end{equation}
the elements $n_p$ will be determined recursively. Put $s^{-1}\left(
w\right) =s^{-1}\cdot w\cdot s,w\in N.$ We have 
\begin{equation}
v\cdot s^{-1}\left( u\right) =\left( \tau \cdot \overrightarrow{\prod }%
n_p\right) \cdot v_0\cdot s^{-1}\left( \overleftarrow{\prod }n_p^{-1}\right)
.  \label{nnn}
\end{equation}
We shall say that an element $x\in {\cal G}$ is in the big cell in ${\cal G}$
if, for all values of the argument $z,$ the value $x\left( z\right) $ is in
the big Bruhat cell $B\bar N\subset G.$

\begin{lemma}
$v\cdot s^{-1}\left( u\right) $ is in the big cell in ${\cal G}$ and admits
a factorization 
\[
v\cdot s^{-1}\left( u\right) =x_{+}^1\cdot x_{-}^1,\;x_{+}^1\in {\cal N}%
,\;x_{-}^1\in \overline{{\cal N}}
\]
\end{lemma}

Indeed, let $u=u_{h/2}u_{h/2-1}\cdot \cdot \cdot u_1,u_p\in {\cal N}^p,$ be
a similar decomposition of $u;$ then we have simply $x_{-}=s^{-1}\left(
u_1\right) .\ $(It is clear that $x_{+}^1\in {\cal B}$ actually does not
have an ${\cal H}$-component and so belongs to ${\cal N}.)$

A comparison of the r.h.s in (\ref{nnn}) with the Bruhat decomposition of
the l.h.s. immediately yields that the first factor in (\ref{n}) is given by 
$n_1=s\left( x_{-}\right) ^{-1}.$

Assume that $n_1,n_2,...,n_{k-1}$ are already computed. Put 
\[
m_k=n_1n_2\cdot \cdot \cdot n_{k-1} 
\]
and consider the element 
\begin{equation}
L^k:=s^{-k+1}\left( \tau \left( m_k^{-1}\right) \cdot \left( v\cdot
s^{-1}\left( u\right) \right) \cdot s^{-1}\left( m_k\right) \right) .
\label{Lk}
\end{equation}

\begin{lemma}
$L^k$ is in the big cell in ${\cal G}$ and admits a factorization 
\begin{equation}
L^k=x_{+}^kx_{-}^k,\;x_{+}^k\in {\cal N},\;x_{-}^k\in \overline{{\cal N}}.
\label{fk}
\end{equation}
\end{lemma}

The elements $x_{\pm }^k$ are computed recursively from the known
quantities. By applying a similar transform to the r.h.s. of (\ref{nnn}) we
get 
\begin{eqnarray}
L^k &=&s^{-k+1}\left( \tau \left( m_k^{-1}\right) \cdot \left( \tau \cdot 
\overrightarrow{\prod_p}n_p\right) \cdot v_0\cdot s^{-1}\left( 
\overleftarrow{\prod_p}n_p^{-1}\right) \cdot s^{-1}\left( m_k\right) \right)
=  \label{Lrhs} \\
&&  \nonumber
\end{eqnarray}
\[
s^{-k+1}\left( \left( \tau \cdot \overrightarrow{\prod_{p\geq k}}n_p\right)
v_0\right) \cdot s^{-k}\left( \overleftarrow{\prod_{p\geq k+1}}%
n_p^{-1}\right) \cdot s^{-k}\left( n_k^{-1}\right) . 
\]
A comparison of (\ref{Lrhs}) and (\ref{fk}) yields $x_{-}^k=s^{-k}\left(
n_k^{-1}\right) ;$ hence $n_k=s^k\left( x_{-}^k\right) ,$ which concludes
the induction.

\section{Generalized Miura Transform}

Our construction of the Miura transform for q-difference operators may be
regarded as a nonlinear version of the corresponding construction for
differential operators, due to Drinfeld and Sokolov \cite{DS}. Recall that
the space of abstract differential operators associated with a given
semisimple Lie algebra ${\frak g}$ is realized as the quotient space of the
affine manifold ${\cal M}^f=f+L{\frak b}\subset L{\frak g}$ (the translate
of $L{\frak b}$ by a fixed principal nilpotent element $f\in {\frak n}$)
over the gauge action of $LN.$ The cross-section theorem of Drinfeld and
Sokolov provides a global model ${\cal S}$ for this quotient space. It is
easy to see that the affine submanifold $f+L{\frak h}\subset f+L{\frak b}$
is a {\em local} cross-section of the gauge action $LN\times $ ${\cal M}%
^f\rightarrow {\cal M}^f$ (i.e., the orbits of $LN$ are transversal to $f+L%
{\frak h})$ and hence $f+L{\frak h}$ provides a local model of the same
quotient space. Thus we get a Poisson structure on $f+L{\frak h}$ and a
Poisson mapping $f+L{\frak h}\rightarrow {\cal S}.$ The computation of the
induced Poisson structure on $f+L{\frak h}$ follows the general prescription
of Dirac (see, e.g., \cite{Flato}), but is in fact greatly simplified, since
all correction terms in the Dirac formula identically vanish. One may notice
that the affine manifold $f+L{\frak h}\subset {\cal M}^f$ is the
intersection of the level surfaces of two moment maps associated with the
gauge actions of the opposite triangular subgroups $LN$ and $L\overline{N};$
it is this symmetry between $LN$ and $L\overline{N}$ that accounts for
cancellations in the Dirac formula. The situation in the nonlinear case is
exactly similar.

We pass to the formal description of our construction. Let $\overline{B}%
\subset G$ be the opposite Borel subgroup, $\overline{N}\subset \overline{B}$
its nilradical, $\overline{{\cal B}}=L\overline{B},\overline{{\cal N}}=L%
\overline{N}.$ Let us consider the Poisson reduction of the space ${\cal C}$
of q-difference connections over the action of the opposite gauge group $%
\overline{{\cal N}}.$ We equip ${\cal C}$ with the Poisson structure (\ref
{taubrack}), where the choice of $\theta $ may be arbitrary.

\begin{proposition}
(i) The q-gauge action $\overline{{\cal N}}\times {\cal C}\rightarrow {\cal C%
}$ leaves $\overline{{\cal B}}\subset {\cal C}$ invariant. (ii) Hamiltonian
vector fields on ${\cal C}$ generated by gauge invariant functions $\varphi
\in C^\infty \left( {\cal C}\right) ^{\overline{{\cal N}}}$ are tangent to $%
\overline{{\cal B}}\subset {\cal C}.$
\end{proposition}

\begin{corollary}
$\overline{{\cal B}}/\overline{{\cal N}}\subset {\cal C}/\overline{{\cal N}}$
is a Poisson submanifold.
\end{corollary}

{\em Remark.} Heuristically, the submanifold $\overline{{\cal B}}\subset 
{\cal C}$ corresponds to reduction at the `zero level' of the moment, hence
the constraints are automatically of the first class.$\square $

We shall now define an embedding $i:{\cal H}\rightarrow {\cal M}^s\cap 
\overline{{\cal B}}$ into the intersection of two 'level surfaces`.

Let $w_0\in W$ be the longest element; let $\pi \in Aut$ $\Delta _{+}$ be
the automorphism defined by $\pi \left( \alpha \right) =-w_0\cdot \alpha
,\alpha \in \Delta _{+}.$ Let $N_i\subset N$ be the 1-parameter subgroup
generated by the root vector $e_{\pi \left( \alpha _i\right) },\alpha _i\in
P.$ Choose an element $u_i\in N_i,u_i\neq 1.$

\begin{lemma}
\label{imm}\cite{st} $w_0u_iw_0^{-1}\in Bs_iB.$
\end{lemma}

We may choose $u_i$ in such a way that $w_0u_iw_0^{-1}\in Ns_iN.$ Put $%
x=u_lu_{l-1}...u_1;$ then $f:=$ $w_0xw_0^{-1}\in Ns^{-1}N\cap \bar N.$

{\em Remark. }The choice of $u_i$ is not unique, however this non-uniqueness
does not affect the arguments below.$\square $

Define the immersion $i:H\rightarrow G:x\longmapsto x\cdot f\cdot s(x^{-1}).$

\begin{proposition}
$i(H)\subset Ns^{-1}N\cap \bar B.$
\end{proposition}

{\em Remark}. By dimension count it is easy to see that $i(H)$ is open in $%
Ns^{-1}N\cap \bar B.$ For $G=SL(n)$ we have simply $i(H)=Ns^{-1}N\cap \bar
B; $ it seems plausible that this is true in the general case as well.$%
\square $

We define the corresponding embedding $i:{\cal H}\rightarrow {\cal G}$ for
loop groups associated with $H,G$ by the same formula. Clearly, $i({\cal H}%
)\subset {\cal M}^s\cap \overline{{\cal B}}.$

\begin{proposition}
$i({\cal H})$ is a local cross-section of the gauge actions ${\cal N}\times 
{\cal M}^s\rightarrow {\cal M}^s,\overline{{\cal N}}\times \overline{{\cal B}%
}\rightarrow \overline{{\cal B}}.$
\end{proposition}

In other words, gauge orbits of ${\cal N},\overline{{\cal N}}$ are
transversal to $i({\cal H})\subset {\cal M}^s\cap \overline{{\cal B}}.$

Let us now assume that the Poisson structure on the space of q-difference
connections is the one described in theorem \ref{main}. We may consider $i(%
{\cal H})$ as a (local) model of the reduced space ${\cal M}^s/{\cal N}$
obtained by 'fixing the gauge` by means of the 'subsidiary condition` $L\in 
\overline{{\cal B}},$ or, alternatively, as a model of $\overline{{\cal B}}/%
\overline{{\cal N}}$ obtained by choosing the subsidiary condition $L\in 
{\cal M}^s.$ The choice of $r_0$ assures that both the 'constraints` and the
'subsidiary conditions` are of the first class. The reduced Poisson
structure on $i({\cal H})$ may be expressed in terms of the Dirac bracket.
As it appears, it is possible to avoid the actual computation of the
'correction terms`. We shall prove the following assertion.

\begin{proposition}
\label{reduced}The quotient Poisson structure on $i({\cal H})$ is given by 
\begin{equation}
\left\{ \varphi ,\psi \right\} _{i({\cal H})}=\left\langle \tilde P_{{\cal H}%
}\nabla \varphi ,\nabla \psi \right\rangle ,\tilde P_{{\cal H}}=\frac{\left(
Id-\tau \right) \left( Id-R_s\right) }{Id-R_s\cdot \tau }.  \label{cartan}
\end{equation}
\end{proposition}

It will be convenient to introduce another parametrization of the Cartan
subgroup which is related to the twisted factorization problem (\ref{twist})
in ${\cal G}.$

\begin{lemma}
$\overline{{\cal B}}\subset {\cal G}^{\prime }.$
\end{lemma}

{\em Proof.} The twisted factorization problem in $\overline{{\cal B}}$ (cf.
(\ref{twist}) amounts to the relation 
\[
\bar b=x_{+}\cdot x_{-}^{-1}n_{-}^{-1},\text{where }x_{+}\in {\cal H}%
,\,x_{-}\in {\cal H},n_{-}\in \overline{{\cal N}}\text{ and }x_{-}=s\left(
x_{+}\right) , 
\]
or, equivalently, 
\begin{equation}
\bar b=x\cdot s(x)^{-1}n_{-}^{-1}.  \label{fAc}
\end{equation}

The same assertion of course holds true for ${\cal H}\subset \overline{{\cal %
B}};$ in that case we have $n_{-}=1.$ Let $\pi :B\rightarrow H$ be the
projection map which assigns to $h\in B$ the element $x\in H$ satisfying (%
\ref{fAc}).

\begin{lemma}
Let $H\subset G$ be the Cartan subgroup. The mapping $p:H\rightarrow
H:x\longmapsto x\cdot s(x)^{-1}$ is an immersion.
\end{lemma}

Put 
\begin{equation}
P_{{\cal H}}=\frac{R_s\cdot \left( \tau -Id\right) }{\left( Id-R_s\right)
\left( Id-R_s\cdot \tau \right) }  \label{lambda}
\end{equation}
and define the Poisson bracket on ${\cal H}$ by 
\begin{equation}
\left\{ \varphi ,\psi \right\} _{{\cal H}}=\left\langle P_{{\cal H}}D\varphi
,D\psi \right\rangle .  \label{final}
\end{equation}

\begin{lemma}
$p:\left( {\cal H},\left\{ ,\right\} _{P_{{\cal H}}}\right) \rightarrow
\left( {\cal H},\left\{ ,\right\} _{i({\cal H})}\right) $ is a Poisson
mapping.
\end{lemma}

Hence to prove proposition \ref{reduced} we may use the Poisson bracket (\ref
{final}) instead of (\ref{cartan}). Let $\varphi ,\psi \in C^\infty ({\cal H}%
).$ Let $\varphi ^{*}=\varphi \circ \pi ,,\psi ^{*}=\psi \circ \pi \in
C^\infty (\overline{{\cal B}})$ be their lifts to $\overline{{\cal B}}$
defined via the twisted factorization map. (In other words, 
\[
\varphi ^{*}\left( \bar b\right) =\varphi \left( x\right) ,\text{ where }%
\bar b=x\cdot s(x)^{-1}n_{-}^{-1},x\in {\cal H}.) 
\]
In the right trivialization of the cotangent bundle of $\overline{{\cal B}}$
the differential $d\varphi ^{*}\left( h\right) \in {\bf b}^{*}$ of $\varphi
^{*}$ is given by 
\begin{equation}
d\varphi ^{*}\left( h\right) =-\tau ^{-1}r_{-}\nabla \varphi ,  \label{dphi}
\end{equation}
where $\nabla \varphi \in {\bf h}$ is the right invariant differential of $%
\varphi $ evaluated at $x=\pi \left( h\right) ,$ and similarly for $\psi
^{*}.$ The standard embedding ${\bf b}^{*}\subset {\bf g}$ allows to regard $%
d\varphi ^{*}\left( h\right) $ as an element of ${\bf g.}$ To compute the
Poisson bracket $\left\{ \varphi ,\psi \right\} $ we may apply proposition 
\ref{reduce}. We have 
\begin{equation}
\left\{ \varphi ,\psi \right\} \left( \pi \left( h\right) \right)
=\left\langle P_{{\cal C}},d\varphi ^{*}\left( h\right) \wedge d\psi
^{*}\left( h\right) \right\rangle .
\end{equation}
Using formula (\ref{Pbracket}) and inserting the expression (\ref{dphi}) for
the differentials we get (\ref{final}).

Let $c:{\cal M}^s\rightarrow {\cal S}$ be the canonical mapping which
assigns to each $L\in {\cal M}^s$ the unique element $L^0\in {\cal S}{\cal \ 
}$lying in the same ${\cal N}$-orbit. {\em The generalized Miura transform }$%
{\bf m}$ is defined by ${\bf m}=c\circ i:{\cal H}\rightarrow {\cal S}.$

\begin{theorem}
The generalized Miura transform is a Poisson mapping.
\end{theorem}

The Poisson structure in ${\cal H}$ is given by (\ref{final}); the Poisson
structure in the target space is the reduced Poisson structure described in
theorem \ref{main}.The proof immediately follows from the fact that $i\left( 
{\cal H}\right) $ and ${\cal S}$ are different models of the quotient space $%
{\cal M}^s/{\cal N}.$

Note that for $G=SL(2)$ our construction of the Miura transform coincides
with the one described in \cite{I}.

\section{The $SL(n)$ case}

Our aim in this section is to compare the Poisson structures arising via the
q-Drinfeld-Sokolov reduction with the results in \cite{FR}. (The case of $n=2$ has been discussed in detail in \cite{I}. Our
analysis for general $n$ is parallel to that of \cite{I}, Section 3,
though our conventions are slightly different.) To begin with,
let us list the standard facts concerning the structure of $SL(n).$ We keep
to the choice of order in the root system of ${\frak sl}(n)$ made in section
1, that is, positive root vectors correspond to lower triangular matrices.
We order simple roots in such a way that the Coxeter element $s=s_1s_2\cdot
\cdot \cdot s_{n-1}$ is acting on the Cartan subalgebra ${\frak h}$ as a
cyclic permutation, 
\[
s_{}^{-1}\cdot diag\left( H_1,H_2,...,H_n\right) =diag\left(
H_n,H_1,...,H_{n-1}\right) ; 
\]
its representative in $G=SL(n)$ is given by 
\[
s^{-1}=\left( 
\begin{array}{ccccc}
0 & -1 & 0 & \cdot \cdot \cdot & 0 \\ 
0 & 0 & -1 & \cdot \cdot \cdot & 0 \\ 
\cdot \cdot \cdot & \cdot \cdot \cdot & \cdot \cdot \cdot & \cdot \cdot \cdot
& \cdot \cdot \cdot \\ 
0 & 0 & 0 & \cdot \cdot \cdot & -1 \\ 
1 & 0 & 0 & \cdot \cdot \cdot & 0
\end{array}
\right) . 
\]
The automorphism $\pi =-w_0$ of ${\frak h}$ is given by 
\[
\pi \cdot diag\left( H_1,H_2,...,H_n\right) =diag\left(
-H_n,-H_{n-1},...,-H_1\right) . 
\]
We may choose the unipotent elements $u_i,$ $i=1,2,...,n-1,$ in such a way
that the principal nilpotent element $f$ constructed in lemma \ref{imm} is
given by 
\[
f=\left( 
\begin{array}{ccccc}
1 & -1 & 0 & \cdot \cdot \cdot & 0 \\ 
0 & 1 & -1 & \cdot \cdot \cdot & 0 \\ 
\cdot \cdot \cdot & \cdot \cdot \cdot & \cdot \cdot \cdot & \cdot \cdot \cdot
& \cdot \cdot \cdot \\ 
0 & 0 & 0 & \cdot \cdot \cdot & -1 \\ 
0 & 0 & 0 & \cdot \cdot \cdot & 1
\end{array}
\right) ; 
\]
the manifold ${\cal M}^s$ consists of matrices of the form 
\[
L=\left( 
\begin{array}{ccccc}
\ast & -1 & 0 & \cdot \cdot \cdot & 0 \\ 
\ast & * & -1 & \cdot \cdot \cdot & 0 \\ 
\cdot \cdot \cdot & \cdot \cdot \cdot & \cdot \cdot \cdot & \cdot \cdot \cdot
& \cdot \cdot \cdot \\ 
\ast & * & * & \cdot \cdot \cdot & -1 \\ 
\ast & * & * & \cdot \cdot \cdot & *
\end{array}
\right) . 
\]
Let $x=diag(x_1,x_2,...,x_n);$ the embedding $i:{\cal H}\rightarrow {\cal M}%
^s$ defined in proposition \ref{imm} is given by 
\[
i(x):=\Lambda =\left( 
\begin{array}{ccccc}
x_1x_n^{-1} & -1 & 0 & \cdot \cdot \cdot & 0 \\ 
0 & x_2x_1^{-1} & -1 & \cdot \cdot \cdot & 0 \\ 
\cdot \cdot \cdot & \cdot \cdot \cdot & \cdot \cdot \cdot & \cdot \cdot \cdot
& \cdot \cdot \cdot \\ 
0 & 0 & 0 & \cdot \cdot \cdot & -1 \\ 
0 & 0 & 0 & \cdot \cdot \cdot & x_nx_{n-1}^{-1}
\end{array}
\right) . 
\]
It is convenient to introduce affine coordinates on $i\left( {\cal H}\right) 
$ in the following way: 
\[
\Lambda (z)=\left( 
\begin{array}{ccccc}
\Lambda _1\left( z\right) \  & -1 & 0 & \cdot \cdot \cdot & 0 \\ 
0 & \Lambda _2\left( qz\right) & -1 & \cdot \cdot \cdot & 0 \\ 
\cdot \cdot \cdot & \cdot \cdot \cdot & \cdot \cdot \cdot & \cdot \cdot \cdot
& \cdot \cdot \cdot \\ 
0 & 0 & 0 & \cdot \cdot \cdot & -1 \\ 
0 & 0 & 0 & \cdot \cdot \cdot & \Lambda _n\left( q^{n-1}z\right) \ 
\end{array}
\right) . 
\]
Let $L_{can}=m(\Lambda )\in {\cal S}$ be the canonical form of $\Lambda ,$%
\[
L_{can}=\left( 
\begin{array}{ccccc}
0 & -1 & 0 & \cdot \cdot \cdot & 0 \\ 
0 & 0 & -1 & \cdot \cdot \cdot & 0 \\ 
\cdot \cdot \cdot & \cdot \cdot \cdot & \cdot \cdot \cdot & \cdot \cdot \cdot
& \cdot \cdot \cdot \\ 
0 & 0 & 0 & \cdot \cdot \cdot & -1 \\ 
1 & u_1\left( z\right) & u_2\left( z\right) & \cdot \cdot \cdot & 
u_{n-1}\left( z\right)
\end{array}
\right) ; 
\]
we put $u_p\left( z\right) :=\left( -1\right) ^{n-p-1}s_{n-p}\left(
q^{n-p}z\right) .$

\begin{proposition}
(\cite{I}, Lemma 2)
We have 
\begin{equation}
s_p\left( z\right) =\sum_{1\leq j_1<j_2<...<j_p\leq n}\Lambda _{j_1}\left(
z\right) \Lambda _{j_2}\left( qz\right) \cdot \cdot \cdot \Lambda
_{j_p}\left( q^{p-1}z\right) .  \label{q-s}
\end{equation}
\end{proposition}

\begin{proposition}
The Poisson bracket on ${\cal H}$ is given by 
\begin{equation}
\left\{ \Lambda _p\left( z\right) ,\Lambda _p\left( w\right) \right\}
=\sum_{m=-\infty }^\infty \frac{\left( 1-q^m\right) \left( 1-q^{m\left(
n-1\right) }\right) }{1-q^{nm}}\left( \frac zw\right) ^m\Lambda _p\left(
z\right) \Lambda _p\left( w\right) ,  \label{FR}
\end{equation}
\[
\left\{ \Lambda _p\left( z\right) ,\Lambda _s\left( w\right) \right\}
=\sum_{n=-\infty }^\infty \frac{\left( 1-q^m\right) \left( 1-q^{-m}\right) }{%
1-q^{nm}}\left( \frac zw\right) ^m\Lambda _p\left( z\right) \Lambda _s\left(
w\right) ,s>p.
\]
\end{proposition}

{\em Proof}. Let $\omega =\exp \frac{2\pi i}n$ be the primitive root of
unity. The eigenvectors of $s$ in ${\frak h}$ are 
\[
e_k=diag(1,\omega ^{-k},...,\omega ^{-(n-1)k}),s\cdot e_k=\omega
^ke_k,k=1,...,n-1. 
\]
The kernel of the Poisson operator $\tilde P_{{\cal H}}=\frac{\left( Id-\tau
\right) \left( Id-R_s\right) }{Id-\tau \cdot R_s}$ is given by the formal
Laurent series 
\[
\sum_{m=-\infty }^\infty \sum_{k=1}^{n-1}\frac 1n\frac{\left( 1-q^m\right)
\left( 1-\omega ^k\right) }{1-q^n\omega ^k}\left( \frac zw\right)
^me_k\otimes e_{n-k}. 
\]
We have 
\begin{eqnarray}
&&\left\{ \Lambda _p\left( z\right) ,\Lambda _s\left( w\right) \right\}
\label{bra} \\
&=&\sum_{m=-\infty }^\infty \sum_{k=1}^{n-1}\frac 1n\frac{\left(
1-q^m\right) \left( 1-\omega ^k\right) }{1-q^m\omega ^k}\left( \frac
zw\right) ^mq^{m\left( s-p\right) }\left( e_k\cdot \Lambda \left( z\right)
\right) _{pp}\cdot \left( e_{n-k}\cdot \Lambda \left( w\right) \right) _{ss}
\nonumber   \\
&=&\sum_{m=-\infty }^\infty \sum_{k=1}^{n-1}\frac 1n\frac{\left(
1-q^m\right) \left( 1-\omega ^k\right) }{1-q^m\omega ^k}\left( \frac
zw\right) ^mq^{m\left( s-p\right) }\omega ^{k\left( s-p\right) }\Lambda
_p\left( z\right) \Lambda _s\left( w\right) .  \nonumber
\end{eqnarray}
Observe that 
\begin{equation}
\frac 1n\sum_k\frac{\left( 1-\omega ^k\right) }{1-q^m\omega ^k}\omega
^{k\left( s-p\right) }=\left\{ 
\begin{array}{c}
\frac{1-q^{n\left( N-1\right) }}{1-q^{mn}},\text{if }s=p, \\ 
q^{-m\left( s-p\right) }\frac{1-q^{-m}}{1-q^{mn}},\text{if }s>p.
\end{array}
\right.  \label{fraction}
\end{equation}
Substituting (\ref{fraction}) into (\ref{bra}), we get (\ref{FR}).

Formula (\ref{q-s}) coincides with the $q$--deformed Miura
transformation defined in \cite{FR}. Formula (\ref{FR}) coincides with
the Poisson bracket on $\Lambda_i(z)$'s derived in
\cite{FR}. Therefore in the case of ${\frak sl}(n)$ the Poisson
algebra obtained by the difference Drinfeld-Sokolov reduction
coincides with the $q$--deformed ${\cal W}$--algebra introduced in
\cite{FR}.

\end{document}